\documentclass[a4paper,12pt]{article}
\usepackage{amsmath, amsfonts, amssymb, amsthm}

\newcounter{num}[section]

\newenvironment{theorem}
{\refstepcounter{num}%
\bigskip\noindent\nopagebreak[4]{\bf Theorem~\arabic{section}.\arabic{num}. }\it}
\newenvironment{lemma}
{\refstepcounter{num}%
\bigskip\noindent\nopagebreak[4]{\bf Lemma~\arabic{section}.\arabic{num}. }\it}
\newenvironment{remark}
{\refstepcounter{num}%
\bigskip\noindent\nopagebreak[4]{\bf Remark~\arabic{section}.\arabic{num}. }}

\begin{document}
\title{Equationally extremal semilattices}
\author{Artem N. Shevlyakov}

\maketitle

\abstract{In the current paper we study extremal semilattices with respect to their equational properties. In the class $\mathbf{S}_n$ of all semilattices of order $n$ we find semilattices which have maximal (minimal) number of consistent equations. Moreover, we find a semilattice in $\mathbf{S}_n$  with maximal sum of numbers of solutions of equations. }

\section{Introduction}

In the current paper we study equations over finite semilattices. Let $\mathbf{S}_n$ denote the class of all semilattices of order $n$. The class $\mathbf{S}_n$ contains extremal semilattices relative to their equational properties. For example, $\mathbf{S}_n$ contains a semilattice $S$ with a maximal (minimal) number of consistent equations. 

Below we study such extremal semilattices from $\mathbf{S}_n$. In Theorem~\ref{th:1} we find a semilattice with a maximal (minimal) number of consistent equations. The semilattice found is isomorphic to the linearly ordered semilattice (to the fan semilattice, respectively). In Theorem~\ref{th:2} we prove that the fan semilattice has a maximal value of 
\[
\Sigma(S)=\sum_{\mathbf{eq}}|{\mathrm{V}}_S(\mathbf{eq})|,
\]
where $\mathbf{eq}$ iterates the set of all equations in one variable. Unfortunately, the corresponding semilattice with a minimal sum is not found yet.  We can only formulate the following conjecture

\medskip

\noindent{\bf Conjecture.} {\it Let $\Sigma(S)$ has the minimum value $\Sigma$ in the set $\mathbf{S}_n$. Then there exists a semilattice $S_1\in \mathbf{S}_n$ with a unique co-atom such that $\Sigma(S_1)=\Sigma$.}

\medskip

We discuss this conjecture in Section~\ref{sec:discuss}

\section{Basic notions}
\label{sec:basics}

Let $S$ be a semilattice of order $n$, where the meet operation is denoted by $\cdot$ and the minimal element is denoted by $0$. The partial order over a semilattice $S$ is defined by $\cdot$ as follows 
\[
x\leq y\Leftrightarrow xy=x.
\]

Moreover, we assume that {\it any semilattice below contains the maximal element $1$}. This assumption is not strict, since one can adjoint $1$ to any semilattice with no maximal element. An element $s\in S$ is called an {\it atom (co-atom)} if there is not an element $t\in S$ with $0<t<s$ (respectively, $s<t<1$). The set of all atoms of a semilattice $S$ is denoted by ${\mathrm{At}}(S)$. Let ${\mathrm{Anc}}(s)$ denote the set of all ancestors of an element $s\in S$, i.e. 
\[
t\in{\mathrm{Anc}}(s) \Leftrightarrow t<s \mbox{ and } \nexists u \; t<u<s.
\]

Let $S$ be a semilattice and $X=\{x_1,x_2,\ldots,x_m\}$ a variable set. An {equation \it over $S$}  is an ordered pair $(t(X)s_1,s(X)s_2)$ or $(t(X)s_1,s_2)$, where  $t(X),s(X)$ are nonempty products of variables of $X$ and $s_i\in S$. Below we shall write equations in traditional form: $t(X)s_1=s(X)s_2$, $t(X)s_1=s_2$. 

\begin{remark}
Let us give some remarks about the set of equations over $S$:
\begin{enumerate}
\item Since an equation is an {\it ordered} pair, the expressions $t(X)s_1=s(X)s_2$,  $s(X)s_2=t(X)s_1$ are assumed to be different equations.

\item Formally, the expression $s_1=t(X)s_2$ is not an equation over $S$, but $t(X)s_2=s_1$ is.

\item Since $1\in S$, the following expressions $t(X)s_1=s(X)$, $t(X)=s(X)s_2$, $t(X)=s(X)$, $t(X)=s_2$ are equations over $S$. 
\end{enumerate}
\end{remark}

An equation of the form $t(X)s_1=s(X)s_2$ ($t(X)s_1=s_2$) is called an equation of the {\it first (respectively, second) kind}. Let ${\mathrm{Eq}}_m^i(S)$ denote the set of all equations of the $i$-th kind in at most $m$ variables over a semilattice $S$.
 
Let ${\mathrm{Eq}}_m(S)={\mathrm{Eq}}_m^1(S)\cup{\mathrm{Eq}}_m^2(S)$ be the set of all equations in at most $m$ variables over a semilattice $S$. For a semilattice $S$ of order $n$ one can directly compute that $|{\mathrm{Eq}}_m^1(S)|=(2^m-1)^2n^2$, $|{\mathrm{Eq}}_m^2(S)|=(2^m-1)n^2$, and finally 
\[
|{\mathrm{Eq}}_m(S)|=(2^m-1)^2n^2+(2^m-1)n^2=2^mn^2(2^m-1).
\]

\section{Two series of semilattices}

Let $L_n=\{a_0,a_1,a_2,\ldots,a_{n-1}\}$ be the linearly ordered semilattice of order $n$ ($n\geq 2$). The order over $L_n$ is given by
\[
a_0<a_1<\ldots<a_{n-1}.
\]

Let $F_n=\{0,a_1,a_2,\ldots,a_{n-2},1\}$ be the semilattice consisting of $n$ ($n\geq 3$) elements, and $A_n=\{a_i\mid 1\leq i\leq n-2\}$ is an anti-chain. 

\begin{center}
\begin{picture}(300,100)

\put(65,80){$F_n\colon$}
\put(100,0){\line(3,1){120}}
\put(100,0){\line(-1,1){40}}
\put(100,0){\line(1,2){20}}
\put(85,0){$0$}

\put(125,40){$a_2$}

\put(225,40){$a_{n-2}$}
\put(160,40){$\ldots$}
\put(65,40){$a_1$}

\put(220,40){\line(-3,1){120}}
\put(60,40){\line(1,1){40}}
\put(120,40){\line(-1,2){20}}

\put(105,80){$1$}
\end{picture} 
\nopagebreak

Fig.1
\end{center}

\section{Inconsistent equations}

In~\cite{shevl_random_semilattices} it was described extremal semilattices of order $n$ with minimal (maximal) number of inconsistent equations. Let ${\mathrm{Eq}}_{m\emptyset}(S)$ be the set of all inconsistent equations over a semilattice $S$ in at most $m$ variables  

\begin{theorem}\textup{\cite[Theorem~1]{shevl_random_semilattices}}
\label{th:1}
For any semilattice $S$ of order $n$ it holds
\begin{equation}
|{\mathrm{Eq}}_{m\emptyset}(L_n)|\leq|{\mathrm{Eq}}_{m\emptyset}(S)|\leq |{\mathrm{Eq}}_{m\emptyset}(F_n)|,
\label{eq:double_ineq_for_Eq_m(emptyset)}
\end{equation}
and  
\[
|{\mathrm{Eq}}_{m\emptyset}(L_n)|=(2^m-1)\frac{n(n-1)}{2},\; |{\mathrm{Eq}}_{m\emptyset}(F_n)|=(2^m-1)(n^2-3n+3)
\]
\end{theorem}

Moreover, in~\cite{shevl_random_semilattices} it was proved that the empty set is most probable in a random generation of equations. More formally, let ${\mathrm{Eq}}_1(S,Y)\subseteq{\mathrm{Eq}}_1(S)$ be the set of all equations with the solution set $Y\subseteq S$, and we have the following statement.

\begin{theorem}\textup{\cite[Theorem~4]{shevl_random_semilattices}}
For any semilattice $S$ of order $n\geq 6$ and any algebraic subset $Y\subseteq S$ it holds
\[
|{\mathrm{Eq}}_1(S,Y)|\leq|{\mathrm{Eq}}_1(S,\emptyset)|.
\]
\end{theorem}

\section{Number of all solutions}

Let us consider an expression
\[
\Sigma(S)=\sum_{\mathbf{eq}\in{\mathrm{Eq}}_1(S)}|{\mathrm{V}}_S(\mathbf{eq})|,
\]
where an equation $\mathbf{eq}$ iterates the set of all equations in one variable,
and find a semilattice with minimal (maximal) value of $\Sigma(S)$. Firstly, one can directly show that
\[
\Sigma(S)=\sum_{s\in S}|{\mathrm{cov}}(s)|,
\]
where ${\mathrm{cov}}(s)={\mathrm{cov}}_1(s)\cup{\mathrm{cov}}_2(s)$ and ${\mathrm{cov}}_i(s)=\{\mathbf{eq}\in{\mathrm{Eq}}_1^i(S)\mid s\in{\mathrm{V}}_S(\mathbf{eq})\}$. 

\begin{lemma}
For any element $s$ of a semilattice $S$ with $|S|=n$ it holds
\begin{equation}
{\mathrm{cov}}_2(s)=\{xa=b\mid a\in S, b=as\},\; |{\mathrm{cov}}_2(s)=n|.
\end{equation}
\label{l:cov_2(s)=n}
\end{lemma}
\noindent {\bf Proof.} An equation $xa=b$ of the second kind is satisfied by $s$ iff $b=as$. Therefore the right part of the equation is determined by a constant $a$ of the left part. Since $a$ belongs to the $n$-element set $S$, there are exactly $n$ equations in $\mathrm{cov}_2(s)$. $\blacktriangleleft$

Thus, $|{\mathrm{cov}}(s)|=|{\mathrm{cov}}_1(s)|+n$ and below we study only the set ${\mathrm{cov}}_1(s)$. Let 
\[
\uparrow s=\{t\mid t\geq s\},\; \perp s=\{t\mid ts=0\}.
\]

\begin{lemma}
\label{l:big}
For any semilattice $S$ of the order $n$ we have
\begin{enumerate}
\item 
\begin{equation}
{\mathrm{cov}}_1(0)=\{xa=xa^\prime\mid a,a^\prime\in S\},\; |{\mathrm{cov}}_1(0)|=n^2;
\label{eq:cov_1(0)}
\end{equation}
\begin{equation}
{\mathrm{cov}}_1(1)=\{xa=xa\mid a\in S\}\; |{\mathrm{cov}}_1(1)|=n;
\label{eq:cov_1(1)}
\end{equation}
\item for any $s\in S$ 
\begin{equation}
{\mathrm{cov}}_1(s)=\bigcap_{t\in {\mathrm{Anc}}(s)}{\mathrm{cov}}_1(t)
\label{eq:cov_1(many_anc)}
\end{equation}
if $|{\mathrm{Anc}}(s)|>1$,
and
\begin{equation}
{\mathrm{cov}}_1(s)={\mathrm{cov}}_1(s^\prime)\setminus \{xa=xa^\prime\mid a\geq s, a^\prime\geq s^\prime, a^\prime\ngeq s\}
\label{eq:cov_1(one_anc)}
\end{equation}
if ${\mathrm{Anc}}(s)=\{s^\prime\}$;

\item for any $s\leq t$ it holds
\begin{equation}
\label{eq:monotonicity}
{\mathrm{cov}}_1(s)\leq{\mathrm{cov}}_1(t)
\end{equation}

\item if $a\in{\mathrm{At}}(S)$ then 
\begin{equation}
{\mathrm{cov}}_1(a)=\{xs=xt\mid s,t\geq a\}\cup\{xs=xt\mid s,t\ngeq a\},
\label{eq:cov_1(atom)}
\end{equation}
\begin{equation}
|{\mathrm{cov}}_1(a)|=|\uparrow a|^2+|\perp a|^2=|\uparrow a|^2+(n-|\uparrow a|)^2.
\label{eq:|cov_1(atom)|}
\end{equation}
\end{enumerate}
\end{lemma}
\noindent {\bf Proof.}
\begin{enumerate}
\item The proof is straightforward. 

\item We prove only~(\ref{eq:cov_1(many_anc)}), since the proof of~(\ref{eq:cov_1(one_anc)}) is straightforward. If $xa=xa^\prime\in\mathrm{cov}_1(s)$ then 
\[
sa=sa^\prime\Rightarrow tsa=tsa^\prime\Rightarrow ta=ta^\prime\Rightarrow xa=xa^\prime\in \mathrm{cov}_1(t)
\]
for any $t\in\mathrm{Anc}(s)$. Thus,
\[
{\mathrm{cov}}_1(s)\subseteq \bigcap_{t\in {\mathrm{Anc}}(s)}{\mathrm{cov}}_1(t).
\]
Let us prove the converse inclusion. Suppose $xa=xa^\prime\in \bigcap_{t\in {\mathrm{Anc}}(a)}{\mathrm{cov}}_1(t)$ but $sa\neq sa^\prime$. 
Then there exist $t_1,t_2\in {\mathrm{Anc}}(s)$ such that $sa\leq t_1$, $sa^\prime\leq t_2$. The equalities $t_ia=t_ia^\prime$ imply 
\[
sa=(sa)t_1=(st_1)a=t_1a=t_1a^\prime,\; sa^\prime=(sa^\prime)t_2=(st_2)a^\prime=t_2a^\prime=t_2a.
\]
Therefore,
\[
aa^\prime s=a^\prime(as)=a^\prime t_1a^\prime=t_1a^\prime=sa,
\]
\[
aa^\prime s=a(a^\prime s)= a t_2 a=t_2 a=sa^\prime,
\]
and we obtain a contradiction $sa=sa^\prime$.

\item Directly follows from the previous statement.

\item The statement immediately follows from formulas~(\ref{eq:cov_1(0)},\ref{eq:cov_1(one_anc)}), since  ${\mathrm{Anc}}(s)=\{0\}$ for any atom $s$.$\blacktriangleleft$
\end{enumerate}

The following lemmas contain results about the sets ${\mathrm{cov}}_1(s)$ in semilattices $L_n,F_n$.

\begin{lemma}
For each $s\in L_n$ we have
\begin{equation}
{\mathrm{cov}}_1(s)=\{xa=xa^\prime\mid a,a^\prime\geq s\}\cup\{xa=xa\mid a<s\},
\label{eq:cov_1(s)_over_L_n}
\end{equation}

\begin{equation}
|{\mathrm{cov}}_1(a_i)|=(n-i)^2+i,
\label{eq:|cov_1(s)|_over_L_n}
\end{equation}

\end{lemma}

\begin{lemma}
For each $s\in F_n$, $n\geq 4$ we have
\begin{multline}
{\mathrm{cov}}_1(a_i)=\{xa_j=xa_k\mid j,k\in \{0,1,2,\ldots,n-2\}\setminus\{i\}\}\\
\cup\{xa_j=xa_k\mid j,k\in\{i,n-1\}\}\; (i\notin \{1,n\}).
\label{eq:cov_1(s)_over_F_n}
\end{multline}
\begin{equation}
|{\mathrm{cov}}_1(a_i)|=(n-2)^2+4.
\label{eq:|cov_1(s)|_over_F_n}
\end{equation}
\end{lemma}

\begin{theorem}
\label{th:2}
For any semilattice $S$ of order $n\geq 4$ we have
\begin{equation}
\sum_{s\in S}|{\mathrm{cov}}_1(s)|\leq \sum_{s\in F_n}|{\mathrm{cov}}_1(s)|.
\label{eqq:rrrrrrr}
\end{equation}
\end{theorem}

\noindent {\bf Proof.} It is clear that for any $s\in S\setminus\{1\}$ we have
\begin{equation}
|\uparrow s|\geq 2.
\label{eqq:yyyyyyyy}
\end{equation}
Suppose $S$ contains at least two atoms then for each atom $a$ we have
\begin{equation}
2\leq |\perp a_i|. 
\label{eqq:yyyyyyyy2}
\end{equation}
and~(\ref{eq:|cov_1(atom)|}) gives
\begin{equation}
|{\mathrm{cov}}_1(a)|\leq (n-2)^2+4
\end{equation} 
for each $a\in{\mathrm{At}}(S)$.

By~(\ref{eq:monotonicity}), for each $s\in S\setminus\{0,1\}$ it holds
\begin{equation}
|{\mathrm{cov}}_1(s)|\leq (n-2)^2+4
\end{equation}
and, using~(\ref{eq:cov_1(0)},\ref{eq:cov_1(1)},\ref{eq:cov_1(s)_over_F_n}), we have
\[
\sum_{s\in S}|{\mathrm{cov}}_1(s)|=|{\mathrm{cov}}_1(0)|+|{\mathrm{cov}}_1(1)|+\sum_{s\in S\setminus\{0,1\}}|{\mathrm{cov}}_1(s)|\leq
n^2+n+(n-2)((n-2)^2+4)=\sum_{s\in F_n}|{\mathrm{cov}}_1(s)|.
\]

\bigskip

Suppose now that $S$ has exactly one atom $a$. Then~(\ref{eq:cov_1(atom)}) gives 
\begin{equation}
|{\mathrm{cov}}_1(a)|=(n-1)^2+1.
\label{eqq:atom_a}
\end{equation}

Let us prove the inequality
\begin{equation}
|{\mathrm{cov}}_1(s)|\leq (n-2)^2+2
\label{eqq:aaaaaaaaaaa}
\end{equation}
for each $s\in S\setminus\{0,a\}$. We have exactly two cases.

\begin{enumerate}
\item Indeed, if $a$ has a unique descendant $a^\prime$ (i.e. $\mathrm{Anc}(a^\prime)=\{a\}$) one can directly compute that 
\[
{\mathrm{cov}}_1(a^\prime)=\{xb=xc\mid b,c\geq a^\prime\}\cup\{xa=xa,x0=x0\},\; |{\mathrm{cov}}_1(a^\prime)|=(n-2)^2+2,
\]
and, using~(\ref{eq:monotonicity}), we obtain~(\ref{eqq:aaaaaaaaaaa}).

\item Suppose $a$ is the ancestor of $a_1,a_2,\ldots,a_k$ ($k\geq 2$). Then one can compute that 
\[
{\mathrm{cov}}_1(a_i)=\{xb=xc\mid b,c\geq a_i\}\cup\{xb=xc\mid b,c\in \perp_a a_i\}\cup\{x0=x0\},
\]
where $\perp_a a_i=\{s\in S\mid sa_i=a\}$. Thus,
\[
|{\mathrm{cov}}_1(a_i)|=|\uparrow a_i|^2+|\perp_a a_i|^2+1.
\]
We have 
\begin{equation}
\label{eqq:uparrow}
2\leq|\uparrow a_i|\leq n-3,
\end{equation}
(since $a_i,1\in\uparrow a_i$ and $0,a,a_j\notin \uparrow a_i$ for $j\neq i$),
\begin{equation}
\label{eqq:perp_a}
2\leq|\perp_a a_i|\leq n-3,
\end{equation}
(since $a,a_j\in\perp_a a_i$ for $j\neq i$, and $0,1,a_i\notin \perp_a a_i$). 
Since $\uparrow a_i\cap\perp_a a_i=\emptyset$ and $0\notin \left(\uparrow a_i\cup\perp_a a_i\right)$ then
\begin{equation}
\label{eqq:leq_n-1}
|\uparrow a_i|+|\perp_a a_i|\leq n-1.
\end{equation}
By~(\ref{eqq:uparrow}--\ref{eqq:leq_n-1}), one can directly compute that 
\[
|\uparrow a_i|^2+|\perp_a a_i|^2\leq (n-3)^2+2^2=n^2-6n+13.
\]
Therefore, 
\[
|{\mathrm{cov}}_1(a_i)|=|\uparrow a_i|^2+|\perp_a a_i|^2+1\leq n^2-6n+14.
\]
For $n\geq 4$ the inequality $n^2-6n+14\leq (n-2)^2+2$ holds, and we obtain~(\ref{eqq:aaaaaaaaaaa}).
\end{enumerate}

Using~(\ref{eq:cov_1(0)},\ref{eq:cov_1(1)},\ref{eqq:atom_a},\ref{eqq:aaaaaaaaaaa}), we obtain
\begin{multline*}
\sum_{s\in S}|{\mathrm{cov}}_1(s)|=|{\mathrm{cov}}_1(0)|+|{\mathrm{cov}}_1(1)|+|{\mathrm{cov}}_1(a)|+
\sum_{s\notin\{0,a,1\}}|{\mathrm{cov}}_1(s)|\leq
n^2+n+(n-1)^2+1+\\
+(n-3)\left((n-2)^2+2\right)=
n^2+n+(n-2)\left((n-2)^2+4\right)=\sum_{s\in F_n}|{\mathrm{cov}}_1(s)|
\end{multline*}
that concludes the theorem. $\blacktriangleleft$

\section{Discussing the conjecture}
\label{sec:discuss}

In Introduction we mentioned the following conjecture.

\medskip

\noindent{\bf Conjecture.} {\it Let $\Sigma(S)$ has the minimum value $\Sigma$ on the set $\mathbf{S}_n$. Then there exists a semilattice $S_1\in \mathbf{S}_n$ with a unique co-atom and $\Sigma(S_1)=\Sigma$.}

\medskip

Surprisingly, the linearly ordered semilattice $L_n$ has not the minimal value $\Sigma(S)$ in $\mathbf{S}_n$. For example, $\Sigma(S)$ has the minimum  at the following semilattice $S_5$ in $\mathbf{S}_5$.

\begin{center}
\begin{picture}(200,120)
\put(100,0){\line(1,1){30}}
\put(100,0){\line(-1,1){30}}
\put(70,30){\line(1,1){30}}
\put(130,30){\line(-1,1){30}}
\put(100,60){\line(0,1){40}}

\put(110,0){25}
\put(80,30){13}
\put(135,30){13}
\put(105,60){7}
\put(105,100){5}
\put(75,100){$S_5$:}
\end{picture} 
\nopagebreak

Fig.2
\end{center} 
(the number at each vertex $s$ is $|\mathrm{cov}_1(s)|$), but the linearly ordered semilattice $L_5$ has the following values of $|\mathrm{cov}_1(s)|$
\begin{center}
\begin{picture}(200,150)
\put(100,0){\line(0,1){120}}
\put(105,0){25}
\put(105,30){17}
\put(105,60){11}
\put(105,90){7}
\put(105,120){5}
\put(75,120){$L_5$:}
\end{picture} 
\nopagebreak

Fig.3
\end{center} 

We have 
\[
\Sigma(L_5)=25+17+11+7+5>25+13+13+7+5=\Sigma(S_5).
\]
Using a computer, we checked the conjecture for each $n\in[4,10]$.

The information of the author:

Artem N. Shevlyakov

Sobolev Institute of Mathematics

644099 Russia, Omsk, Pevtsova st. 13

Phone: +7-3812-23-25-51.

e-mail: \texttt{a\_shevl@mail.ru}
\end{document}